\documentclass[graybox]{svmult}

\usepackage{type1cm}        
\usepackage{makeidx}         
\usepackage{graphicx}        
\usepackage{multicol}        
\usepackage[bottom]{footmisc}

\usepackage{newtxtext}       %
\usepackage[varvw]{newtxmath}       

\makeindex             

\usepackage{siunitx}
\usepackage{amsbsy}
\usepackage{amsmath,bm}
\usepackage{amsfonts}

\usepackage{multirow}
\usepackage{tabu}

\usepackage[bottom]{footmisc}
\usepackage{cite}

\usepackage{url}

\begin{document}

\title*{Numerical Assessment of PML Transmission Conditions in a Domain Decomposition\\ Method for the Helmholtz Equation}
\titlerunning{PML Transmission Conditions for the Helmholtz Equation}
\author{Niall Bootland 
\and
Sahar Borzooei 
\and
Victorita Dolean 
\and
Pierre-Henri Tournier 
}
\authorrunning{N. Bootland et al.}
\institute{Niall Bootland \at STFC Rutherford Appleton Laboratory, Harwell Campus, UK,
	\email{niall.bootland@stfc.ac.uk}
	\and Sahar Borzooei \at University Côte d'Azur, CNRS, LJAD, France,
	\email{Sahar.Borzooei@univ-cotedazur.fr}
	\and Victorita Dolean \at University Côte d'Azur, CNRS, LJAD, France, and University of Strathclyde, UK,
	\email{work@victoritadolean.com}
	\and Pierre-Henri Tournier \at Sorbonne Université, CNRS, Université Paris Cité, Inria, Laboratoire Jacques-Louis Lions (LJLL), F-75005 Paris, France,
	\email{pierre-henri.tournier@sorbonne-universite.fr}}

%
%
\maketitle

\abstract{The convergence rate of domain decomposition methods (DDMs) strongly depends
on the transmission condition at the interfaces between subdomains. Thus, an important aspect in improving the efficiency of such solvers is careful design of appropriate transmission conditions.
In this work, we will develop an efficient solver for Helmholtz equations based on perfectly matched layers (PMLs) as transmission conditions at the interfaces within an optimised restricted additive Schwarz (ORAS) domain decomposition preconditioner, in both two and three dimensional domains. 
We perform a series of numerical simulations on a model problem and will assess the convergence rate and accuracy of our solutions compared to the situation where impedance boundary conditions are used.}


\section{Introduction}
\label{Borzooei:Intro}

Finite element discretisations of large-scale time-harmonic wave problems typically lead
to ill-conditioned linear systems with a large number of unknowns. A promising class of methods to solve such huge systems in parallel, both in terms of convergence and computing time, is offered by domain decomposition methods (DDMs). These approaches rely on a partition of the computational domain into smaller
subdomains, leading to subproblems of smaller sizes which are manageable by direct solvers. A~robust domain decomposition (DD) preconditioner for large-scale computations is given in \cite{Borzooei-b15}.
However, improving the efficiency of such preconditioners continues to be a challenging issue.
Recent work has shown that transmission operators
based on perfectly matched layers (PMLs) are well-suited for two-dimensional configurations of the Helmholtz problem within non-overlapping DDMs \cite{Borzooei-v1}. Further, PMLs have been used successfully as transmission conditions in DDMs applied to geophysical applications modelled by the Helmholtz equation \cite{Borzooei-tournier}.

In this work, we present an efficient PML-based Schwarz-type preconditioner with overlapping subdomains to solve large-scale wave propagation problems. 
We then assess the performance of this one-level DD algorithm, where the transmission conditions at the boundaries between subdomains are PML conditions in order to provide a better approximation to the transparent boundary operator. Further, we will investigate the convergence properties and compare them with the use of more standard impedance transmission conditions.
 
\section{Mathematical model}
\label{Borzooei:sec_2}

As an underlying model we consider the Helmholtz equation in free space, given by
\begin{equation}
	\label{Borzooei:eq_Helm}
	- (\Delta + k^2(\mathbf{x})) u(\mathbf{x}) = g(\mathbf{x}), \quad  \mathbf{x} \in \Omega
\end{equation}
for $\Omega = \mathbb{R}^{d}$ in dimension $d = 2$ or $3$, where $k(\mathbf{x}) = \frac{2\pi}{\lambda}$ is the wavenumber, with $\lambda = \frac{c}{f}$ being the wavelength, $c(\mathbf{x})$ the wave speed and $f$ the frequency. Note that the angular frequency is then defined as $\omega = 2\pi f$. To close the problem we prescribe the physically relevant condition at infinity known as the far field Sommerfeld radiation condition
\begin{align}
	\label{Borzooei:eq_Sommerfeld}
	\lim_{|\mathbf{x}|\rightarrow\infty} |\mathbf{x}|^{\frac{d-1}{2}} \left( \frac{\partial u}{\partial |\mathbf{x}|} - iku \right) = 0.
\end{align}
Since we can not compute on the whole free space domain, we consider truncating to an appropriate finite domain. Let us suppose now that $\Omega \subset \mathbb{R}^{d}$ represents a finite computational domain capturing the physical area of interest. A typical approach, as in \cite{Borzooei-sah}, is to replace the Sommerfeld condition \eqref{Borzooei:eq_Sommerfeld} with the first-order approximation
\begin{align}
	\label{Borzooei:eq_Imp}
	\frac{\partial u}{\partial \mathbf{n}} + iku = 0, \quad \mathbf{x} \in \partial \Omega,
\end{align}
known as the impedance (or Robin) boundary condition (Imp BC), with $\mathbf{n}$ being the unit outward normal to the boundary $\partial\Omega$. This enables the appropriate description of wave behaviour in a bounded domain. The finite element discretisation of \eqref{Borzooei:eq_Helm} can then be written as a linear system $A \mathbf{u} = \mathbf{b}$.

\subsection{PML formulation}
\label{borzooei:subsec_1}

Perfectly matched layers (PMLs) were introduced as a better alternative to absorbing boundary conditions (ABCs) by Berenger\cite{Borzooei-1} to achieve a higher accuracy in domain truncation by eliminating undesired numerical reflections from boundaries, leading to exponential convergence of the numerical solution to the exact solution \cite{Borzooei-4}. PML implementation is done by stretching Cartesian coordinates such that the stretching is defined in a layer surrounding $\Omega$, as in \cite{Borzooei-4a}, giving a larger computational domain~$\Omega_\text{PML}$. In this regard, we assume the boundaries of the artificially truncated domain~$\Omega$ are aligned with the coordinate axes.

For simplicity of exposition, we will focus on truncating the problem in the $x$ direction. Let us suppose that the PML extends from the boundary of our domain of interest at $x=a$ to $x = a^{*}$ and  Dirichlet conditions are imposed on $x = a^{*}$. The stretched coordinate mapping used is given by
\begin{align}
	\label{Borzooei:eq_pml}
	\frac{\partial}{\partial x_\text{pml}} \mapsto \frac{1}{1 - \frac{i}{\omega}\sigma(x)} \frac{\partial }{ \partial x}, \quad \text{where} \quad
	\begin{cases}
		\begin{array}{ll}
			\sigma(x) = 0 & \text{ if } x < a, \\
			\sigma(x) > 0 & \text{ if } a < x < a^{*}.
		\end{array}
	\end{cases}
\end{align}

In the PML region, where $\sigma(x)>0$, oscillating solutions turn into exponentially
decaying ones. In the original domain $\Omega$, $\sigma(x) = 0$ so that the underlying equation is unchanged. In this work we will study three different stretching functions \cite{Borzooei-4aa, Borzooei-4ac}, namely
\begin{align}
	\label{Borzooei:eq_sigmas}
	\sigma_{-1}(x) = \frac{1}{a^{*}-x}, \quad \sigma_{-2}(x) = \frac{2}{(a^{*}-x)^{2}}, \quad \text{and} \quad \sigma_{2}(x) &= \alpha(a^{*}-x)^{2}.
\end{align}
In $\sigma_{2}(x)$, $\alpha$ is experimentally chosen to take the value 30 in our simulations. To incorporate a PML into other coordinate directions, we simply apply equivalent one-dimensional transformations to obtain $\frac{\partial}{\partial y_\text{pml}}$ and $\frac{\partial}{\partial z_\text{pml}}$. At the corners of the extended computational domain $\Omega_\text{PML}$ we will have PML regions that stretch along two or three directions simultaneously; this will not generate any problems.
Implementing this mapping in, for instance, a three-dimensional domain requires a slight change to the Helmholtz equation~\eqref{Borzooei:eq_Helm} over $\Omega_\text{PML}$, resulting in the Laplace operator $\Delta$ being replaced by following operator which stretches in the PMLs
\begin{align}
	\label{Borzooei:eq_lap}
	\Delta_\text{pml} = \frac{\partial^{2}}{{\partial x_\text{pml}}^2} + \frac{\partial^{2}}{{\partial y_\text{pml}}^2} + \frac{\partial^{2}}{{\partial z_\text{pml}}^2}.
\end{align}

\subsubsection{Accuracy assessment for PMLs}

\begin{figure}[t!]
	\centering
	\includegraphics[width=0.4\textwidth]{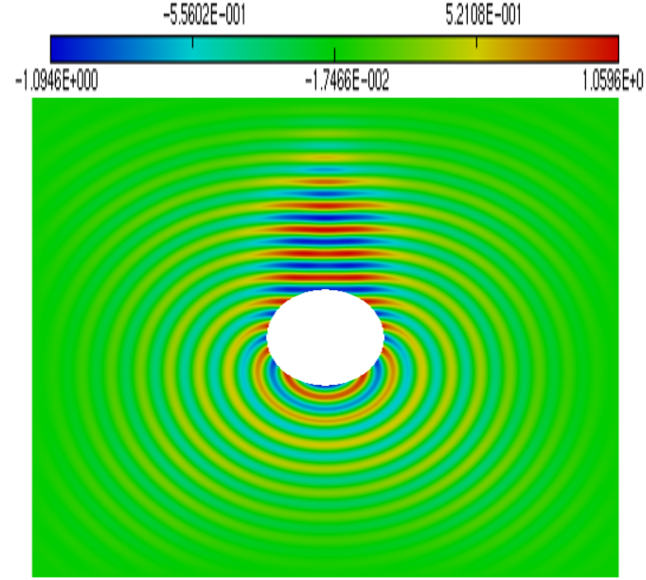}
	\caption{Plane wave excitation solution when using PMLs as global boundary conditions in 2D.}
	\label{Borzooei:fig_PMLglob}   
\end{figure}

In this section we will solve the Helmholtz equation when PMLs are applied as global boundary conditions for a 2D domain of length $10\lambda$ in each direction. We will compute the $L^2$ relative error with respect to the analytical exact solution and compare it with the situation where impedance boundary conditions are used instead. We consider a scattering problem of a plane wave by a circular obstacle, with a Dirichlet boundary condition on the boundary of the obstacle, shown in Figure~\ref{Borzooei:fig_PMLglob}. First, in Table~\ref{Borzooei:tab_globalerror}, we compare different stretching functions $\sigma$ with the utilization of higher order P3 Lagrange finite elements and discretization of $n_{\lambda} = 20$ points per wavelength. We find that the best accuracy is obtained with $\sigma_{-1}$ and so we continue our tests with this function here.
Within our tests we vary the number of points per wavelength $n_{\lambda}$ and the PML length in order to investigate their relative effect on the resulting error. Results are detailed in Table~\ref{Borzooei:tab_global}. We see that, except for $n_{\lambda} = 5$, PMLs provide higher accuracy compared to impedance boundary conditions, even when the length of the PMLs incorporate only $0.1\lambda$. Moreover, for a fixed PML length, and again even for $0.1\lambda$, the error still decreases when increasing $n_{\lambda}$, whereas for impedance boundary conditions the error is dominated by the domain truncation even for $n_{\lambda} = 5$. Of course, the error also decreases significantly with increasing PML length, all the way down to $3 \times 10^{-5}$ for $n_{\lambda} = 50$ and $10\lambda$.

\begin{table}[t!]
	\caption{$L^2$ relative error for different stretching functions $\sigma$ with PML length $L_{\text{pml}} = \lambda$. The radius of the circular obstacle is $R = \lambda$.}
	\label{Borzooei:tab_globalerror}
	\centering
	\tabulinesep=0.7mm
	\begin{tabu}{ccc}
		\multicolumn{3}{c}{Stretching functions}\\
		$\sigma_{-1}$ & $\sigma_{-2}$ & $\sigma_{2}$ \\
		\hline
		0.00112  &0.001517  & 0.075495 \\
	\end{tabu}
\end{table}

\begin{table}[t!]
	\caption{$L^2$ relative error for different PML lengths with $\sigma_{-1}$ or impedance boundary conditions (Imp BCs), $R= \lambda$.}
	\label{Borzooei:tab_global}
	
	\smallskip
	
	\centering
	\scalebox{.93}{		
 \tabulinesep=0.7mm
	\begin{tabu}{c|cccccccccc|c}
		& \multicolumn{10}{c|}{PML length} & \\
	  $n_{\lambda}$ &$0.1\lambda$& $0.2\lambda$& $0.3\lambda$ & $0.5\lambda$ &$\lambda$ & $2\lambda$ & $3\lambda$ & $4\lambda$ &$5\lambda$ & $10\lambda$& Imp BCs\\
		\hline
 5&0.10408 &0.02441 &0.02684 &0.02120  & 0.01685& 0.01251&0.00927 &0.00741  & 0.00605&0.00265 &0.05118 \\
 10& 0.01354 & 0.01011& 0.00665&  0.00534&0.00425 &0.00311 & 0.00235& 0.00184 &0.00150 &0.00067 &0.04642 \\
		  
 20&  0.00893&0.00467 &0.00268 & 0.00159 &0.00112 & 0.00078&0.00059 & 0.00046&0.00038 &0.00017 & 0.04620\\
 30 & 0.00797 &0.00320 & 0.00175&  0.00083& 0.00050& 0.00035&0.00026 & 0.00021 &0.00017 &0.00007 & 0.04620\\ 
 40 &0.00617  &0.00246 &0.00121&0.00056 &0.00029  & 0.00020&0.00015 & 0.00012& 0.00009 &0.00006 & 0.046212\\
 50 & 0.00578 & 0.00192&0.00096 & 0.00041 &0.00020 & 0.00013&0.00010 & 0.00008 & 0.00006&0.00003 &0.046216 \\ 
	\end{tabu}
}	
	\vspace{-.3cm}
	
\end{table}

\subsection{Domain decomposition preconditioner}

A preconditioner $M^{-1}$ is a linear operator whose use aims to reduce ill-conditioning and allow faster convergence of an iterative solver. Usually (but not always) this approximates $A^{-1}$ and has a matrix--vector product that is much cheaper to compute than solving the original linear system. To this end, we employ right preconditioning within GMRES to solve our discretised linear system, namely by solving
\begin{align}
	\label{Borzooei:right}
	AM^{-1} \mathbf{y} = \mathbf{f}, \quad \text{where} \quad  \mathbf{u} &= M^{-1}\mathbf{y}.
\end{align}
Right preconditioning benefits from minimising a residual that is independent of the preconditioner, unlike left-preconditioned GMRES.
Overlapping Schwarz methods come with the advantages of better convergence and easier implementation compared to substructuring methods. Furthermore, contrary to non-overlapping methods, corners do not need specific treatment. Overlapping methods are also a natural choice to consider when using PML transmission conditions, as the added PML can be naturally included in the overlap region. In this work we use the optimised restricted additive Schwarz (ORAS) domain decomposition preconditioner, given~by
\begin{align}
	\label{Borooei:oras}
	M_\text{ORAS}^{-1}= \sum_{s=1}^{N_\text{sub}} R_s^T D_s A_s^{-1} R_s,
\end{align}
where $N_\text{sub}$ is the number of overlapping subdomains $\Omega_s$ into which the domain $\Omega$ is decomposed. To define the matrices present in \eqref{Borooei:oras}, let $\mathcal{N}$ be an ordered set of the unknowns of the whole domain and let $\mathcal{N} = \bigcup_{s=1}^{N_\text{sub}} \mathcal{N}_s$ be its decomposition into the non-disjoint ordered subsets corresponding to the different overlapping subdomains~$\Omega_s$. Further, define $N = \vert \mathcal{N} \vert$ and $N_s = \vert \mathcal{N}_s \vert$. The $N_s \times N_s$ matrices $A_s$ stem from the discretisation of local boundary value problems on $\Omega_s$ with transmission conditions chosen as either Robin
or PML conditions to be implemented at the subdomain interfaces. The $N_s \times N$ matrix $R_s$ is the Boolean restriction matrix from~$\Omega$ to subdomain $\Omega_s$ while $R_s^T$ is then the extension matrix from subdomain $\Omega_s$ to~$\Omega$. The $N_s \times N_s$ diagonal matrices $D_s$ provide a discrete partition of unity, i.e., are such that $\sum_{s=1}^{N_\text{sub}} R_s^T D_s R_s = I$. See, e.g., \cite{Borzooei-6,Borzooei-v2} for further details on such methods. PMLs are introduced as transmission conditions on the interface boundaries of the local subdomains in \cite{Borzooei-v1}. In this approach, the PML region is included strictly inside the overlap, the PML region being the outermost layers within each overlapping domain. This ensures that there is enough overlap for the approach to be efficient and sufficient length of the PML for a good approximation of the interface transmission condition.

\section{Numerical results}

\subsection{PML as transmission conditions for a 2D domain}
\begin{figure}[t!]
	\centering
	\ \hfill
	\includegraphics[height=0.385\textwidth,clip,trim=0mm 0mm 0mm 0mm]{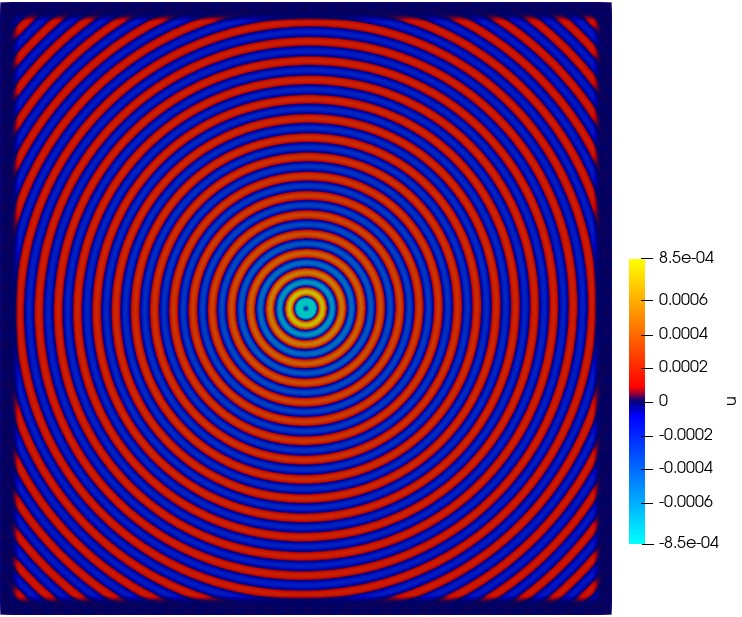}
	\hfill
	\includegraphics[height=0.385\textwidth,clip,trim=0mm 0mm 0mm 0mm]{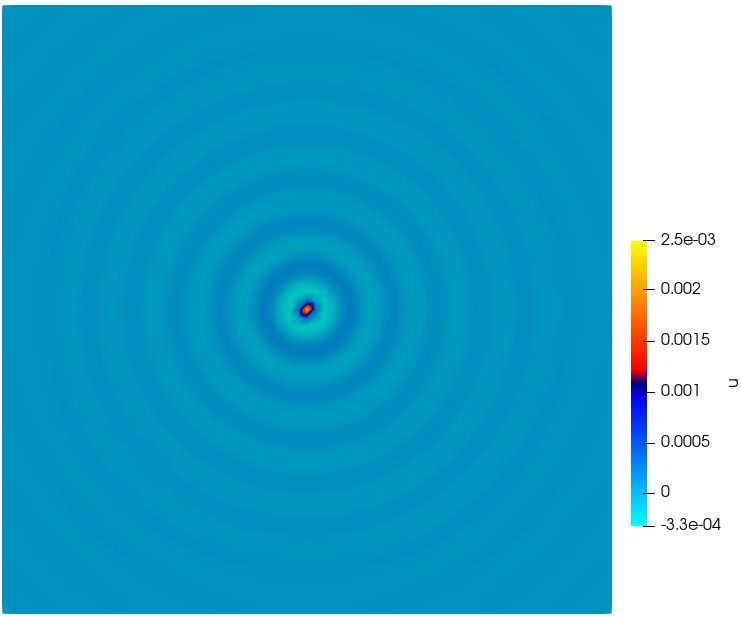}
	\hfill \
	\caption{The solution of the 2D,  with $f=3Hz$ (left), and 3D,  with $f=1Hz$ (right), point source excitation problem with PML boundary conditions.}
	\label{Borzooei:fig_1}
\end{figure}

As a simple model, we consider excitation by a Gaussian point source, $ S(x,y) = e^{-30k((x-5)^2 + (y-5)^2)} $, in the center of a 2D domain of size $[0,10]\times [0,10]$, as shown in Figure~\ref{Borzooei:fig_1} (left). The convergence rate is studied when either PML or impedance conditions are imposed as global boundary conditions (BCs) or interface conditions~(ICs). This leads to four different configurations in total\footnote{In 2D we present results only with PMLs for the global BCs; a comparison with impedance BCs will be given later for the full 3D problem in Table~\ref{Borzooei:tab_3d_1h}.}. To discretise we 
employ~P3 finite elements on regular grids with $n_{\lambda} = 15$.

In our tests, we set the wave speed to $c=1$ and we vary the frequency $f$ from $3$Hz to $10$Hz, which leads to a number of wavelengths in the domain ranging from~30 to~100. That also results in different values of $\#$DoFs, which represents the number of degrees of freedom in the discrete problem. We decompose the global domain into either $N = 8\times8 = 64$ or $N = 10\times10 =100$ square subdomains and use interface PML regions of length $L_\text{pmli}$. 

In Table~\ref{Borzooei:tab_2d_1h}, simulations results are given using the $\sigma_{-1}$ stretching function\footnote{A comparison with other choices of stretching function $\sigma$ will be given later in Table~\ref{Borzooei:tab_3d_vary_sigma}.}, where the interface PML length is $L_\text{pmli} = 1h$ and $h = \frac{\lambda}{n_{\lambda}}$ is the mesh size. The PML length on the global boundary is chosen to be $L_\text{pml} = 2\lambda$ and the number of overlapping layers of elements between subdomains is varied from 2 to 8 layers. We first observe that when interface PMLs are used we always require fewer iterations compared to using impedance ICs. Secondly, with the impedance condition the iteration counts increase with frequency $f$, but this is not the case when using interface PMLs where iteration counts remain insensitive to $f$. Finally, we note that an overlap of 4 layers is sufficient here for the PMLs with little benefit seen as we increase the overlap further while for the impedance condition a larger overlap is needed to continually reduce the iteration counts.

\begin{table}[h!]
	\centering
	\caption{Iteration counts for varying frequency $f$ and choices of ICs, discretised using P3 elements with $n_{\lambda} = 15$. Within the PMLs we use $\sigma_{-1}$, $L_\text{pmli} = 1h$ and $L_\text{pml} = \lambda$.}
	
	\medskip
	
	\tabulinesep=0.7mm
	\begin{tabu}{cc|cc|cccc|cccc}
		& & & & \multicolumn{8}{c}{Overlap} \\
		& & & & \multicolumn{4}{c|}{$N =64$} & \multicolumn{4}{c}{$N = 100$} \\
		BCs & ICs & $f$ (Hz) & $\#$DoFs & 2 & 4 & 6 & 8 & 2 & 4 & 6 & 8\\
		\hline
		PML & Imp & \multirow{2}{*}{3} & \multirow{2}{*}{$2,\!076,\!481$} & 52 & 48  & 44 & 41 & 65 & 58 & 53 & 51 \\
		PML & PML &  & & 39 & 33 & 33 & 32 & 49 & 42 & 40 & 42\\
		\hline
		PML & Imp & \multirow{2}{*}{5} & \multirow{2}{*}{$5,\!480,\!281$} & 63 & 59  & 56  & 54 & 70 & 65 & 60 & 57\\
		PML & PML & & & 41 & 35  & 33 & 34 & 49 & 42 & 41 & 42 \\
		\hline
		PML & Imp & \multirow{2}{*}{10} & \multirow{2}{*}{$21,\!077,\!281$} & 69 & 67 & 63 & 61 & 77 & 71 & 67 & 68 \\
		PML & PML & & & 40 & 34  & 33 & 34 & 48 & 41 & 41 & 41
	\end{tabu}
	\label{Borzooei:tab_2d_1h}  
\end{table}

\begin{table}[h!]


	\centering
	\caption{Iteration counts for $f = 3~Hz$ and varying choices of ICs, discretised using P3 elements with $n_{\lambda} = 15$. Within the PMLs we use $\sigma_{-1}$, $L_\text{pmli} = 5h$ and $L_\text{pml} = \lambda$.}
	
		\medskip
		
	\tabulinesep=0.7mm
	\begin{tabu}{cc|ccc|cccc|cccc}
		& & && & \multicolumn{8}{c}{Overlap} \\
		& & && & \multicolumn{4}{c|}{$N = 64$} & \multicolumn{4}{c}{$N = 100$} \\
		BCs & ICs & $f$ (Hz)&$\#$DoFs & &2 & 4 & 6 & 8 & 2 & 4 & 6 & 8\\
		\hline
		PML & Imp & \multirow{2}{*}{3}& \multirow{2}{*}{$2,\!076,\!481$} && 52 & 48 & 44 & 41 & 65 & 58 & 53 & 51 \\
		PML & PML && & & 110 & 64 & 31 & 31 & 140 & 80 & 39 & 39
	\end{tabu}
	\label{Borzooei:tab_2d_5h}  
\end{table}
In Table~\ref{Borzooei:tab_2d_5h}, the simulations with $f = 3\,$Hz are repeated but now with $L_\text{pmli} = 5h$. Note that for the appropriate transmission of data between subdomains, we should consider the length of the overlap to be larger than  the length of the PML region. This can be seen in Table~\ref{Borzooei:tab_2d_5h} where an overlap of more than 5 layers is required for good convergence. Comparing the number of iterations, when the overlap is sufficient, with those in Table~\ref{Borzooei:tab_2d_1h}, we can see a small improvement in the convergence when using a larger interface PML region. Note that the one-level preconditioner is by nature not robust, in the sense that the number of iterations usually depends on the number of subdomains. That is to say, the number of iterations does not depend only on the quality of the approximation of the absorbing interface conditions, which as we can see from Table~\ref{Borzooei:tab_global} is already good when using PMLs of small length.

\subsection{PML as transmission conditions for a 3D domain}

In this section, 
we consider a similar Gaussian point source excitation in the center of the 3D domain, $ S(x,y,z) = e^{-30k((x-5)^2 + (y-5)^2 + (z-5)^2)} $; see Figure~\ref{Borzooei:fig_1} (right). For this problem we discretise with P2 finite elements and use $n_{\lambda} = 5$ and $L_\text{pml}=2\lambda$. When recording iteration counts in this section, the use of $-$ means the simulation failed due to memory limitations while $\bullet$ indicates a lack of convergence in~2000~iterations. In Table~\ref{Borzooei:tab_3d_1h} we use $\sigma_{-1}$ and compare all four combinations of BCs and ICs when $L_\text{pmli} = 1h$ and $f = 1\,$Hz, this results in \#DoFs = 2,803,221. We observe that using PML rather than impedance conditions reduces iteration counts both when used as BCs or ICs, in particular, when swapping from impedance for both~BCs and~ICs to~PMLs we see at least a $2/3$ reduction in iterations. Furthermore, using PML BCs again provides a somewhat more accurate solution when comparing~$L^{2}$ relative error with respect to the analytical exact solution.  Here, we consider $N = 6\times6\times5 = 180$ and $N = 7\times7\times6 = 294$ subdomains.

\begin{table}[t!]
	\centering
	\caption{Iteration counts and $L^{2}$ error for $f = 1\,$Hz and varying choices of BCs and ICs, discretised using P2 elements with $n_{\lambda} = 5$. Within the PMLs we use $\sigma_{-1}$, $L_\text{pmli} = 1h$ and $L_\text{pml} = 2\lambda$.}
	
		\smallskip
		
	\tabulinesep=0.7mm
	\begin{tabu}{cc|cccc|cccc|c}
		& & \multicolumn{8}{c|}{Overlap} & \\
		& & \multicolumn{4}{c|}{$N = 180$} & \multicolumn{4}{c|}{$N = 294$} & \\
		BCs & ICs & 2 & 4 & 6 & 8 & 2 & 4 & 6 & 8 & $L^2$ relative error \\
		\hline
		Imp & Imp & 40 & 31 & 27 & $-$ & 45 & 35 & 32 & 34 & \multirow{2}{*}{0.211853}  \\
		Imp & PML & 36 & 29 & 26 & $-$ & 42 & 35 & 29 & 31  \\
		PML & Imp & 30 & 22 & 20 & $-$ & 33 & 25 & 23 & 21& \multirow{2}{*}{0.0709828 } \\
		PML & PML & 24 & 20 & 18 & $-$ & 27 & 23 & 20 & 19 
	\end{tabu}
	\label{Borzooei:tab_3d_1h}  
\end{table}

In Table~\ref{Borzooei:tab_3d_vary_Lpmli}, simulations for the full PML case are repeated for different lengths of~$L_\text{pmli}$. Again we see the overlap should be larger than the interface PML length and, when so, iteration counts slowly decrease as $L_\text{pmli}$ increases.

\begin{table}[t!]
	\centering
	\caption{Iteration counts for $f = 1\,$Hz with PML BCs and ICs varying the PML interface length $L_\text{pmli}$, discretised using P2 elements with $n_{\lambda} = 5$. Within the PMLs we use $\sigma_{-1}$ and $L_\text{pml} = 2\lambda$.}
	
		\smallskip
		
	\begin{tabu}{cc|c|cccc|cccc}
		& & & \multicolumn{8}{c}{Overlap} \\
		& & & \multicolumn{4}{c|}{$N = 180$} & \multicolumn{4}{c}{$N = 294$} \\
		BCs & ICs & $L_\text{pmli}$ & 2 & 4 & 6 & 8 & 2 & 4 & 6 & 8 \\
		\hline
		PML & PML & $1h$ & 24 & 20 & 18 & $-$ & 27 & 23 & 20 & 19 \\
		PML & PML & $2h$ & 30 & 19 & 17 & $-$ & 34 & 22 & 19 & 18 \\
		PML & PML & $4h$ & 37 & 21 & 15 & $-$ & 42 & 24 & 17 & 15 \\
		PML & PML & $6h$ & 38 & 21 & 18 & $-$ & 43 & 24 & 21 & 15
	\end{tabu}
	\label{Borzooei:tab_3d_vary_Lpmli}  
\end{table}

Finally, we compare different stretching functions $\sigma$ for the case 
of $L_\text{pmli} = \nobreak 4h$. The results in Table~\ref{Borzooei:tab_3d_vary_sigma} show that\, the best convergence is provided\, when we choose~$\sigma_{-1}$.  While the iteration counts when using $\sigma_{-2}$ have only a small increase, it is always more effective to use $\sigma_{-1}$. The convergence observed for $\sigma_{2}$ is much poorer, demonstrating the importance of choosing a suitable stretching function in order to be advantageous in the domain decomposition preconditioner. In our tests~$\sigma_{-1}$ provided the best choice and justifies its use in our previous simulations.

\begin{table}[t!]
	\centering
	\caption{Iteration counts for $f = 1\,$Hz and varying choice of ICs and PML stretching function $\sigma$, discretised using P2 elements with $n_{\lambda} = 5$. Within the PMLs we use $L_\text{pmli} = 4h$ and {$L_\text{pml} = 2\lambda$}.}
	
		\smallskip
		
	\begin{tabular}{cc|c|cccc|cccc}
		& & & \multicolumn{8}{c}{Overlap} \\
		& & Stretching & \multicolumn{4}{c|}{$N = 180$} & \multicolumn{4}{c}{$N = 294$} \\
		BCs & ICs & function & 2 & 4 & 6 & 8 & 2 & 4 & 6 & 8 \\
		\hline
		PML & Imp & \multirow{2}{*}{$\sigma_{-1}$} & 30 & 22 & 20 & $-$ & 33 & 25 & 23 & 21 \\
		PML & PML & & 37 & 21 & 15 & $-$ & 42 & 24 & 17 & 15 \\
		\hline
		PML & Imp & \multirow{2}{*}{$\sigma_{-2}$} & 33 & 28 & 24 & $-$ & 38 & 30 & 27 & 26 \\
		PML & PML & & $\bullet$ & 33 & 19 & $-$ & $\bullet$ & 38 & 23 & 19 \\
		\hline
		PML & Imp & \multirow{2}{*}{$\sigma_{2}$} & $\bullet$ & $\bullet$ &973  & $-$ & $\bullet$ & $\bullet$ & 1984 & 1201 \\
		PML & PML & & $\bullet$ & $\bullet$ & $\bullet$  & $-$ & $\bullet$ & $\bullet$ & $\bullet$ & $\bullet$
	\end{tabular}
	\label{Borzooei:tab_3d_vary_sigma}
\end{table}

\section{Conclusion} In this work, we have introduced the use of PMLs as interface conditions within an overlapping domain decomposition solver for Helmholtz equations. With the choice of PMLs as interface conditions, better convergence is achieved compared to using impedance conditions. Results on 2D and 3D model problems show the utility of the approach with a suitable choice of stretching function. 


\section*{Acknowledgement}
This work has received funding from the European
Union’s Horizon 2020 under Marie Curie grant agreement No 847581 (COFUND BoostUrCAreer). The authors are grateful to the OPAL infrastructure from Université Côte d’Azur and the Université Côte d’Azur’s Center for High-Performance Computing for providing resources and support.



\end{document}